\documentclass[reprint, amsmath,amssymb, aps]{revtex4-2}
\pdfoutput=1
\usepackage{graphicx}
\usepackage{dcolumn}
\usepackage{bm}

\usepackage{amsthm}
\usepackage{amssymb}
\usepackage{amsfonts}
\usepackage{mathtools}
\usepackage{tikz-cd}
\usepackage{chngcntr}

\usepackage{amsmath}
\usepackage{tikz-cd}
\usepackage{booktabs} 
\usepackage{caption} 
\usepackage{subcaption} 
\usepackage{graphicx}
\usepackage{pgfplots}
\usepackage{tikz}
\usepackage{algpseudocode,algorithm,algorithmicx}
\usepackage{hyperref}
\usepackage{amssymb}
\usepackage{amsfonts}

\definecolor{blue}{HTML}{1F77B4}
\definecolor{orange}{HTML}{FF7F0E}
\definecolor{green}{HTML}{2CA02C}
\newcommand{\Z}{\mathbb{Z}}
\newcommand{\R}{\mathbb{R}}

\begin{document}

\preprint{APS/123-QED}

\title{Point Fields of Last Passage Percolation and Coalescing Fractional Brownian Motions}

\author{Konstantin Khanin}
\author{Liying Li}
\author{Zhanghan Yin}
\affiliation{ Department of Mathematics,  University of Toronto }

\date{\today}

\begin{abstract}
We consider large-scale point fields which naturally appear in the context of the Kardar-Parisi-Zhang (KPZ) phenomenon. Such point fields are geometrical objects formed by points of mass concentration, and by
shocks separating the sources of these points. We introduce similarly defined point fields for processes of coalescing fractional Brownian motions (cfBM). The case of the Hurst index 2/3 is of particular interest for us since, in this case, the power law of the density decay is the same as that in the KPZ phenomenon. In this paper, we present strong numerical evidence that statistical properties of points fields in these two different settings are very similar. We also discuss theoretical arguments in support of the conjecture that they are exactly the same in the large-time limit. This would indicate that two objects may, in fact, belong to the same universality class.
\end{abstract}

\maketitle

\section{Introduction}
The KPZ equation:
\begin{eqnarray}
  \label{eq:KPZ}
  \partial_t h + (\partial_x h)^2  = \partial_{xx} h + F,\\
  F = \text{ space-time white noise}, \nonumber
\end{eqnarray}
describes the motion of growing surfaces that is subject to smoothing effects, slope-dependent growth speed and space-time uncorrelated noise.
In the seminal paper by Kardar, Parisi and Zhang \cite{PhysRevLett.56.889}, it was predicted that the fluctuations of the height function, $h(t,x)$, are of the order of $t^{1/3}$ and the spatial correlation occurs at the scale of $t^{2/3}$.
The $1:2:3$ scaling, known as the KPZ scaling, also arises in many other models including random matrices, random growth models, interacting particle systems, optimal paths/directed polymers in random
environments, randomly forced Burgers equation/Hamilton--Jacobi equations \cite{PhysRevA.16.732,BDJDistributionLengthLongest1999, JohShapeFluctuationsRandom2000a,PSScaleInvariancePNG2002, CatorGroen,
  10.1214/EJP.v11-366, BKSMicroscopicConcavityFluctuation2012,ACQProbabilityDistributionFree2011a, MQRKPZFixedPoint2020, DOVDirectedLandscape2018, BKGlobalSolutionsRandom2018}.

Besides the height function, there is also a geometrical approach to understand the KPZ scaling through the geometrical properties of optimal paths or equivalent objects in these models. Such geometrical objects already arises when
representing the solution to~(\ref{eq:KPZ}) via the Feynman--Kac formula after applying the Hopf--Cole transform $h(t,x) = -2\nu \ln \phi(t,x)$:
\begin{eqnarray}
\label{eq:FK}
\phi(t,x) &=& \int e^{-\frac{1}{2\nu} \Big[ h(t,\gamma_0)  + \int_0^t F(s,\gamma_{s})\,ds  \Big] } W^{t,x}(d\gamma),\\
W^{t,x}(\cdot) &=& \text{ Wiener measure with endpoint $(t,x)$.} \nonumber
\end{eqnarray}

The Gibbs measure on paths~$$P^{t,x}_0(d\gamma) = Z^{-1}e^{-\frac{1}{2\nu} \Big[ h(t,\gamma_0)  + \int_0^t F(s,\gamma_{s})\,ds  \Big] } W^{t,x}(d\gamma)$$ is a \textit{polymer measure} in the random environment given by~$F$.
Although~$P^{t,x}_0$ is random, there exists a deterministic number~$\chi \in [0,1]$, called the \textit{transversal exponent},  such that the probability
$$P_0^{t,x} \big(\max_{0\le s\le t}|\gamma_s-\gamma_t| = O(t^\chi) \big)$$
is close to $1$ for typical environment.
The KPZ scaling corresponds to~$\chi=2/3$; for comparison, in the absence of randomness, i.e., $F=0$, $P^{t,x}_0$ is equivalent to the Wiener measure and thus~$\chi=1/2$.
Naturally, the transversal exponent which describes the large-scale property should not feel the roughness of the environment~$F$; it is believed that a sufficient condition for the $2/3$ transversal exponent is
rapid space-time decorrelation of the random environment.  In the sequel we will assume $F$ to be smooth which describes the large-scale properties with rapid space-time decorrelation.

In the zero temperature limit~$\nu\to 0$, the Gibbs measures will concentrate on geodesics that have a fixed endpoint~$\gamma_t=x$ and minimize the action
\begin{equation}\label{eq:action}
h(0, \gamma_0) + \int_0^t [L(\dot{\gamma}_s) + F(s,x+\gamma_{s})] \, ds, \quad L(p)=\frac{p^2}{2}.
\end{equation}
The Lagrangian~$L$ can be other convex functions, and this optimization problem is the one that occurs in the Lax--Oleinik variational principle that gives the viscous solution to the inviscid Hamilton--Jacobi equation 
\begin{equation*}
\partial_t u +  H(\partial_x u) = F(t,x),
\end{equation*}
where~$H$ is the Legendre dual of $L$.
Compared to polymer measures, the geometry of geodesics is easier to describe because of fewer layer of randomness.
The transversal exponent $\chi$ of the geodesics can be defined in a similar way: let~$\gamma=\gamma^{t,x}$ be the geodesic, then in typical
environment 
$$\max_{0 \le s \le t} |\gamma_s-\gamma_t| = O(t^{\chi}).$$
In general, the models of finding optimal paths in random environments are called \textit{first/last-passage percolation} (FPP/LPP).

In most of the FPP/LPP models, geodesics cannot intersect except at the endpoints.  In the context of Hamilton--Jacobi equations this means the following: let~$\gamma^{1,2}:[0,t]\to\R$ be two
geodesics of~(\ref{eq:action}) (i.e., any perturbation of $\gamma^{1,2}$ will have higher action); then~$\gamma^1(s) = \gamma^2(s)$ can only happen for~$s=0$ or~$t$.  This is due to the convexity of the Lagrangian~$L$. 

The non-intersecting property gives a monotone structure to the geodesics.
In particular, the map~$x \mapsto \gamma^x(0)$ is monotone (non-decreasing), where $\gamma^x$ is the minimizing path of~(\ref{eq:action}) with~$\gamma(t) = x$.
Although~$\gamma^x$ may not be unique, the monotone map is well-defined since the discontinuity points are at most countable.
Such points with more than one minimizers correspond to the formation of shocks.

We can obtain a consistent family of monotone maps~$(\phi^{s,t})_{s<t}$ that satisfies $\phi^{r,s}\circ\phi^{s,t} = \phi^{r,t}
$ for all~$r<s<t$, if we look at infinite geodesics.  More specifically,
fix a large negative $T$ and let~$\gamma^{t,x}$ be the minimizing path of~(\ref{eq:action}) that starts from time~$T$ and terminates at~$(t,x)$.
The family of monotone maps is given by~$\phi^{s,t}_T(x) =
\gamma^{t,x}(s)$, $T<s<t$.
By the principle of dynamic programming, these monotone maps are consistent:
$$\phi_T^{r,s}\circ\phi_T^{s,t} = \phi_T^{r,t}, \quad T \le s<t<r. $$
Sending~$T\to-\infty$ we get rid of the dependence on~$T$, and can think of obtaining these monotone maps from infinite geodesics.
These monotone maps depend on the random environment, and since the environment is space-time stationary, so are the monotone maps; the temporal stationarity means that~$\phi^{s+r,t+r}$ has the same
statistics as~$\phi^{s,t}$ for all~$r$, and the spatial stationarity means that $x\mapsto \phi^{s,t}(x)-x$ is a stationary process.
Since the transversal exponent for infinite geodesics should be the same as the finite ones, we can also see the KPZ scaling in terms of the monotone maps: $|\phi^{-t,0}(x) - x| = O(t^{\chi})$, $\chi=2/3$. 

We are interested in understanding to what extent the monotonicity property and correlation structures determine the value of $\chi$.
More precisely, let~$\phi^{s,t}:\R\to\R$, $s<t$, be a consistent, stationary flow of random monotone maps.
Then, is it true that for
some~$\chi\ge 0$, $|\phi^{-t,0}(x)-x| = O(t^{\chi})$?
How does $\chi$ depend on the distribution of the monotone maps?
Moreover, when~$\chi$ exists, what is the scaling limit of the renormalized monotone maps $[\mathrm{R}_{L,\chi} \phi]^{s,t}(x) =  L^{-\chi}\phi^{Ls,Lt}(L^{\chi}x)$, as $L \to \infty$?
Is this limit uniquely determined by~$\chi$?
We are particularly interested in the case~$\chi=2/3$, since this value of $\chi$ corresponds to the KPZ universality.

Another special case is $\chi = 1/2$.
This case can be studied rigorously since the monotone maps~$\phi^{s,t}$ are independent in time, and the scaling limit is given by the \textit{coalescing Brownian motion} (cBM).
However, from the geometrical perspective described above, neither the~$2/3$ nor $1/2$ exponents should be special;
it should be possible to obtain scaling limits for other values of $\chi$ by varying the temporal decay of correlation of the monotone maps.

Any scaling limit of the renormalization operator $\mathrm{R}_{L,\chi}$ produces its fixed point.
The fixed point for~$\chi=1/2$ is given by the flow of cBM, constructed as follows.
Particles start from every position on the line at time~$s$ and perform independent Brownian motions until collision. When two particles collide, they merge into a new particle which
continues to perform Brownian motion independent of other particles.
For $s<t$, let~$\phi^{s,t}(x)$ be the time-$t$ position of
the particle coming from location~$x$ at time~$s$.
The coalescing construction ensures that~$(\phi^{s,t})_{s<t}$ is a family of monotone maps.
Moreover, due to the memoryless effect (Markov property) of Brownian motions, if we follow the trajectory of one particle, $t \mapsto \phi^{s,t}(x)$, the trajectory is a Brownian path despite
collisions taking place along the way, and it follows from the diffusive scaling of Brownian motions that $|\phi^{0,t}(x) - x| = O(t^{1/2})$ for the flow of cBM.
By the invariance of Brownian motion it is not hard to see that it is a fixed point for~$R_{L,1/2}$. 
Note the different time directions for these models: the flow of cBM is forward in time, while the infinite geodesics from Hamilton--Jacobi equations are backward in time.

The flow of cBM was first rigorously constructed by Arratia \cite{ArratiaCBM}. The most technical point was to show the ``coalescence from infinity'' property, that is, at every
time~$t>0$, there are only countably many particles left at discrete positions.
As a consequence, all the maps~$\phi^{s,t}$ are piecewise constant functions that can be characterized by two discrete point fields
\begin{equation*}
    \ldots < a_{-1}<a_0 < a_1 < \ldots, \quad \ldots < b_{-1} < b_0 < b_1 < \ldots
\end{equation*}
such that~$\phi^{s,t}((a_n,a_{n+1})) = b_n$ for $n \in \Z$. 

\begin{figure}[h]
    \centering
    \includegraphics[width=0.5\textwidth]{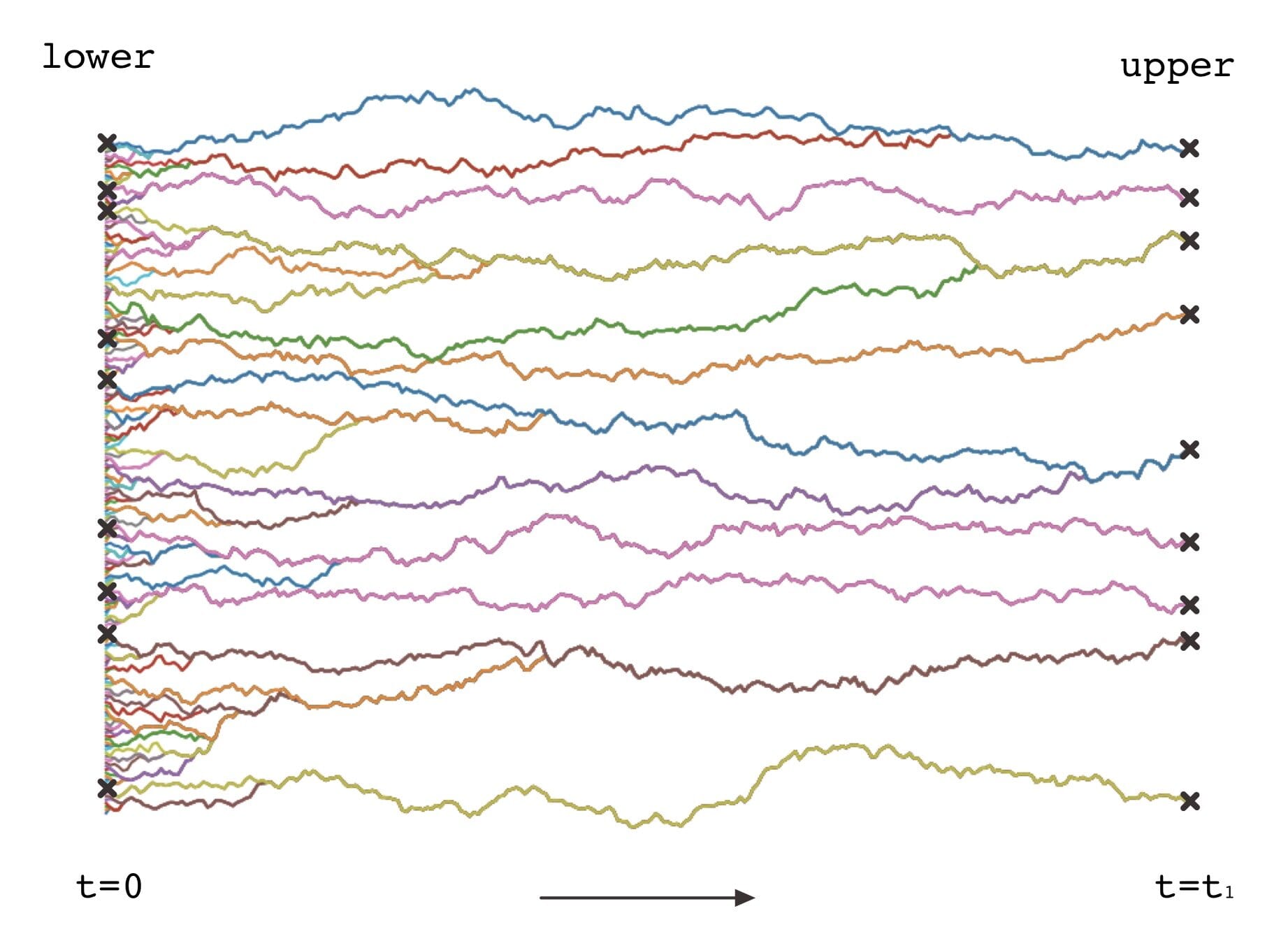}
    \caption{Upper/lower point fields of a coalescing process}
    \label{UnL}
\end{figure}

The ``coalescence from infinity" property means that the random set of surviving particles at any positive time $t>0$ constitutes a point field. We call this the \textit{upper point field}. On the
other hand, for each point in the upper point field, the set of starting positions that end up at that point is almost surely an interval. The endpoints of these intervals constitute the \textit{lower
  point field} (FIG.\ \ref{UnL}).

It is known that for many FPP/LPP models the coalescence of infinite geodesics also takes place \cite{CNSurfaceViewFirstPassage1995, 10.1214/19-AIHP1016}, that is, denoting by~$\gamma^x$ the geodesic from~$x$ at time~$0$, for~$x \neq y$, there is a
(negative) time~$T=T_{x,y}$ for which~$|\gamma^x(T)-\gamma^y(T)| \ll 1$, and $|\gamma^x(t)-\gamma^y(t)|$ will converge exponentially fast for~$t<T$.  In fact, for lattice model, $\gamma^x(t) =
\gamma^y(t) $ for $t < T$.
The time~$T_{x,y}$ is called the coalescence time, and according to the KPZ scaling, for fixed $x$ and~$S$, the starting point~$y$ such that the coalescence time~$|T_{x,y}| < S$ should be
distance~$O(S^{2/3})$ away from~$x$.
As a result, the fixed point of~$R_{L,2/3}$ obtained from solvable KPZ models is also given by piecewise constant maps.
In this paper we argue that the statistical properties of this maps are determined by the monotonicity properties and the planar geometry.

In what follows, we numerically construct various candidates for the fixed point of $\mathrm{R}_{L,2/3}$ from a new class of coalescing processes called the \textit{coalescing fractional Brownian motion} (cfBM). We then compare the statistics of the upper/lower point fields to that obtained from solvable KPZ models, and observe that these models share strikingly similar statistics.

\section{Numerical Experiment}
\subsection{Coalescing Fractional Brownian Motion and Exponential Corner Growth Model}
Intuitively, the  construction of cfBM is similar to that of cBM:
at the initial time $t=0$, independent fractional Brownian particles of Hurst index $H \in (0,1)$ start at every point on $\mathbb{R}$,  and two particles ``coalesce'' into one upon 
 collision.
However, now the dynamics of particles after coalescence admit different interpretations.
 We propose three types of coalescence rules, namely:
\begin{enumerate}
    \item \textbf{Coin-flip: } When two particles collide, one is chosen with equal probability to ``absorb'' the other particle and continue its motion.
    \item \textbf{Regenerate: } When two particles collide, they both vanish and a new independent fBM particle is spawned at the point of collision.
    \item \textbf{P\'olya-urn: } Let $\alpha  \ge 0$ be a fixed \textit{P\'olya index} (reminiscent of P\'olya urn). Every particle starts with weight $1$.
       When two particles of weights $w_1$ and $w_2$ collide, they respectively have probability
      $\frac{w_1^\alpha}{w_1^\alpha + w_2^\alpha}$ and $\frac{w_2^\alpha}{w_1^\alpha + w_2^\alpha}$ of winning.
The winning particle absorbs the losing particle and continues its motion, with a new weight $w_1+w_2$. Note that when $\alpha =
      \infty$, the  particle with higher weight always wins, and when $\alpha=0$, this is the coin-flip model.
\end{enumerate}
For cBM, all the above coalescing rules are equivalent due to the strong Markov property. On the contrary, fBM with Hurst index $ H\neq 1/2$ is non-Markovian \cite[Theorem
2.3]{NouSelectedAspectsFractional2012}, so  one may expect that different coalescing rules would lead to different kinds of dependence on the past, and hence different versions of cfBM. This intuition is supported numerically, see section \ref{Aux}.

Although a rigorous construction of cfBM is currently not available, we study numerical simulations of cfBM with finitely many initial starting points that are sufficiently dense and equally spaced. The details of
the simulation can be found in \cite{SM}.  We are particularly interested in cfBM with $H = 2/3$, where the point fields will have the same density scaling as KPZ models.

We will compare the point fields generated by the cfBM and by the \textit{exponential corner growth model}.
This is a last-passage percolation model known to belong in the KPZ universality class \cite{10.1214/EJP.v11-366}.
In this model, weights of independent and identically distributed exponential random variables
are placed on each $(\Z_+)^2$-lattice point, and a boundary condition is specified on the non-negative $x$- and $y$-axis.  The geodesics are the up-right paths maximizing the sum of weights they visit. For details of this exactly solvable model, see \cite{10.1214/19-AIHP1016}.
The geodesics of this model can be generated efficiently.   For the rest of this paper, by ``LPP'' we refer specifically to this exactly solvable model.

\subsection{Test Statistics}
We are  interested in the following questions:
\begin{enumerate}
\item \textbf{KPZ-like Properties: } Which model of cfBM has upper and/or lower point fields with similar statistics as that of LPP (KPZ)?
\item \textbf{Symmetry:} The duality of cBM says that there is a joint realization of two cBMs --- one forward in time and the other backward, with non-crossing paths. It follows that the upper
      and lower point fields of cBM are identically distributed. The symmetry of upper and lower point fields is also known to hold for LPP \cite{PimDualityCoalescenceTimes2016}.  What about cfBM?   
\item \textbf{Different Coalescing Rules: } Does numerical evidence corroborate with the expectation that different coalescing rules of the cfBM would lead to a difference in  point field statistics?
\end{enumerate}

\begin{figure}[h]
    \centering
    \includegraphics[width=0.5\textwidth]{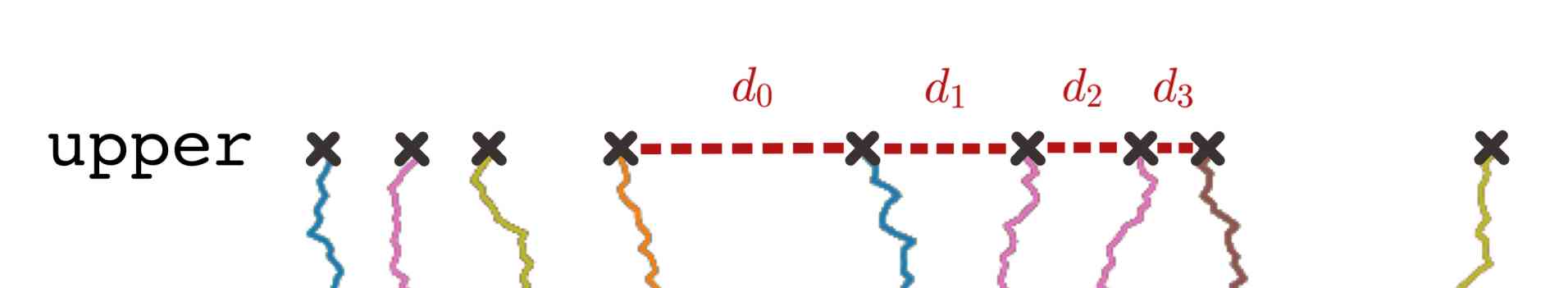}
    \caption{$d_k$ are identically distributed.  $r_k = d_k/d_0$ is the \textit{jump-$k$ ratio}.}
    \label{quants}
  \end{figure}

This involves comparing the upper and lower point fields of different coalescing processes. With the translational invariance, these point fields are characterized by the distributions of distances
between consecutive points and all finite joint-distributions of such distances.
In this paper we will use the $p$-value from the \textit{Kolmogorov-Smirnov} (K-S) test to compare the distributions of the following one-dimensional statistics:
\begin{itemize}
\item \textbf{Normalized Consecutive Point Distance:} We study the distribution of the distance between consecutive points in the point field, normalized by the sample mean.
  We denote this by $\delta_0$.
\item \textbf{Jump-$k$ ratio:} For fixed $k \ge 1$, we consider two intervals of $k$ intervals apart, and call the ratio between their length the \textit{jump-$k$ ratio}, denoted by $r_k$ (FIG.\ \ref{quants}). 
\end{itemize}

The random variables $\delta_0$ and $r_k$, $k \ge 1$, can be defined for  both the upper and the lower point fields. 

\subsection{Results} \label{Aux}
We simulated cfBM with the coalescing rule of coin-flip, regenerate and P\'olya-urn with a few values of $\alpha$.  For the P\'olya-urn model we will present here the data of $\alpha=1$ and
$\alpha=\infty$.  In the supplemental material \cite{SM} a wider range of $\alpha$ are presented.

\begin{table}[H]
  \caption{\label{tab:LPP-vs-cfBM} $p$-values comparing LPP and cfBM models.  See the full table at \cite[TABLE IV]{SM}}
  \centering
\begin{tabular}{|c|c|c|c|c|}
  \hline
  & coin-flip &  $\alpha=1$ & $\alpha=\infty$ & regenerate \\  \hline
  $\delta_0$ & 0.50 & 0.23  & 0.15 & 0.77 \\  \hline
  $r_1$ & 0.86 & 0.45  & 0.59 & 0.45\\  \hline 
  $r_3$ & 0.41 & 0.33  & 0.79 & 0.97\\  \hline 
  $r_6$ & 0.41 & 0.32  & 0.71 & 0.78\\  \hline
\end{tabular}
\end{table}

In TABLE \ref{tab:LPP-vs-cfBM} we compare the upper point fields from cfBM models to that from LPP.
Surprisingly, across the table we obtained relatively high $p$-values, considering that the K-S test is based on the $L^\infty$ distance between 
empirical CDFs that should be sensitive to small differences.
This provides numerical evidence for a strong similarity between the statistical properties of the upper point fields in the KPZ problem and in processes of coalescing fractional Brownian motions with
the Hurst index $2/3$.

We also compute the $p$-values for the lower point fields in \cite[TABLE V]{SM}, and except for the coin-flip model, all the $p$-values are small.
This is expected since as we will see immediately below, coin-flip cfBM and LPP are the
only two models that have symmetry between the upper and lower point fields.
\begin{table}[H]
  \caption{\label{tab:upper-vs-lower} $p$-values of $\delta_0$, $r_k$ between upper and lower fields of every process.  See the full table at \cite[TABLE III]{SM}.}
  \centering
\begin{tabular}{|c|c|c|c|c|c|c|}
    \hline
  & coin-flip &  $\alpha=1$ &  $\alpha=\infty$ & regenerate & LPP\\    \hline
  $\delta_0$ & 0.75  & 0.16 &  $<$0.01 & 0.02 & 0.92\\    \hline
  $r_1$ & 0.51 & 0.67  & $<$0.01 & 0.22 & 0.66\\    \hline
  $r_3$ & 0.55 & 0.51 & $<$0.01 & 0.52 & 0.99\\    \hline
  $r_6 $ & 0.63  & 0.42 & 0.01 & 0.11 & 0.42\\    \hline
\end{tabular}
\end{table}

In TABLE \ref{tab:upper-vs-lower}, we compare the statistics between the upper and lower point fields of each version of cfBM and the LPP model.
As a reference, the upper and lower point fields of the LPP model are known to be identically distributed, and the $p$-values are in the range from $0.42$ to $0.92$.  Using this range as a reference,
TABLE \ref{tab:upper-vs-lower} suggests that the
coin-flip model is the only other model that can also possibly exhibit such symmetry.

There are intuitive reasons why the P\'olya-urn model ($\alpha>0$) and the regenerate models do not have such symmetry in the reverse time direction.
In the P\'olya-urn model, the weights of particles
are monotonely increasing, and in the  regenerate model, the history of particles is erased upon collision.

\begin{table}[H]
\caption{\label{tab:distance}$p$-values of $\delta_0$. (The upper and lower parts of the table correspond to the upper and lower point fields.) See the full table at \cite[TABLE I]{SM}.}\centering
\begin{tabular}{|c|c|c|c|c|}
    \hline
  & coin-flip &  $\alpha=1$&  $\alpha =\infty$ & regenerate  \\    \hline
  coin-flip & &   0.45 &     0.33 &   0.33 \\      \hline
  $\alpha=1$ &  0.17 & &    0.84 &   0.13\\    \hline
    $\alpha=\infty$ &  $<$0.01 &  0.04 & &  0.12 \\    \hline
  regenerate &  $<$0.01 &   $<$0.01 &   $<$0.01 & \\    \hline
\end{tabular}
\end{table}
Lastly, we compare the point fields from cfBM with different coalescing rules.

In TABLE~\ref{tab:distance}, we obtain extremely low $p$-values from the lower point fields. This is expected and confirms that different coalescing rules yield different versions of cfBM.
For the upper point fields, the $p$-values are larger ($>0.10$), which agrees with the data in TABLE \ref{tab:LPP-vs-cfBM}.

In \cite[TABLE II]{SM}, we collect the $p$-values comparing the distribution of the jump-$2$ ratio $r_2$, and similar conclusions can be drawn.

\section{Conclusion}
While most of the research in the field of KPZ was concentrated on the statistics of interface, or so-called, height function, in this paper we suggest a more geometrical point of view. For many
problems in the KPZ universality class one can naturally define the large scale point fields which can be viewed as places of concentration of polymers, or particles, and points of separation which
can be viewed as shocks. The asymptotic statistical properties of these point fields are supposed to be universal, and in some sense encode  the statistical properties of the interface.

In this paper we introduce a new model which was not studied previously --- coalescing fractional Brownian motions. This model is not related to any external disorder setting, and hence, normally would not be considered for belonging to the KPZ universality class. However, in the submitted paper we present theoretical arguments and numerical evidence suggesting that the point fields arising from the model of cfBM have the same statistical properties as  the points fields in the models of the KPZ class. Besides providing a new approach to the problem of KPZ universality, our results suggest that the KPZ universality class is, in fact, much larger than it was previously thought.

It was suggested in \cite{BKGlobalSolutionsRandom2018} that the statistical properties of the point fields are completely determined by the monotonicity properties, decorrelation conditions, and the requirement of a fixed asymptotic power law decay of the density of points fields. The arguments in favour of the above conjecture were based on the renormalization approach. The main idea was that, in the large-time limit, the probability law of a point field converges to a renormalization fixed point which is stable apart from one neutral direction corresponding to different values of the exponent of the power decay of the density. In the present paper we provide a strong numerical support for such a universality of the point fields. \\

We stress again that cfBM can be considered for different Hurst indices. The case of the Hurst index 1/2 corresponds to standard (non-fractional) Brownian motions. Universality in this case was studied in \cite{piterbarg1997expansions,PPIntermittencyTracerGradient1999a,khanin2021asymptotic}. It was shown in \cite{khanin2021asymptotic} that the renormalization procedure can be viewed as the dynamical process of iteration by random monotone piecewise constant maps. The Hurst index $1/2$ corresponds to the situation when maps are identically distributed and independent. It was rigorously proved that in this case the fixed point is stable. 
Other Hurst indices correspond to the situation when different random monotone maps are correlated in time. The rigorous mathematical analysis in this case is a challenging problem.



\providecommand{\noopsort}[1]{}\providecommand{\singleletter}[1]{#1}%
\begin{thebibliography}{3}%
\makeatletter
\providecommand \@ifxundefined [1]{%
 \@ifx{#1\undefined}
}%
\providecommand \@ifnum [1]{%
 \ifnum #1\expandafter \@firstoftwo
 \else \expandafter \@secondoftwo
 \fi
}%
\providecommand \@ifx [1]{%
 \ifx #1\expandafter \@firstoftwo
 \else \expandafter \@secondoftwo
 \fi
}%
\providecommand \natexlab [1]{#1}%
\providecommand \enquote  [1]{``#1''}%
\providecommand \bibnamefont  [1]{#1}%
\providecommand \bibfnamefont [1]{#1}%
\providecommand \citenamefont [1]{#1}%
\providecommand \href@noop [0]{\@secondoftwo}%
\providecommand \href [0]{\begingroup \@sanitize@url \@href}%
\providecommand \@href[1]{\@@startlink{#1}\@@href}%
\providecommand \@@href[1]{\endgroup#1\@@endlink}%
\providecommand \@sanitize@url [0]{\catcode `\\12\catcode `\$12\catcode
  `\&12\catcode `\#12\catcode `\^12\catcode `\_12\catcode `\%12\relax}%
\providecommand \@@startlink[1]{}%
\providecommand \@@endlink[0]{}%
\providecommand \url  [0]{\begingroup\@sanitize@url \@url }%
\providecommand \@url [1]{\endgroup\@href {#1}{\urlprefix }}%
\providecommand \urlprefix  [0]{URL }%
\providecommand \Eprint [0]{\href }%
\providecommand \doibase [0]{https://doi.org/}%
\providecommand \selectlanguage [0]{\@gobble}%
\providecommand \bibinfo  [0]{\@secondoftwo}%
\providecommand \bibfield  [0]{\@secondoftwo}%
\providecommand \translation [1]{[#1]}%
\providecommand \BibitemOpen [0]{}%
\providecommand \bibitemStop [0]{}%
\providecommand \bibitemNoStop [0]{.\EOS\space}%
\providecommand \EOS [0]{\spacefactor3000\relax}%
\providecommand \BibitemShut  [1]{\csname bibitem#1\endcsname}%
\let\auto@bib@innerbib\@empty
\bibitem [{Note1()}]{Note1}%
  \BibitemOpen
  \bibinfo {note} {\protect \url
  {https://github.com/neo-nice-to-meet-you/cfBm-KPZ/}}\BibitemShut {NoStop}%
\bibitem [{\citenamefont {Flynn}(line)}]{PyPI_fbm}%
  \BibitemOpen
  \bibfield  {author} {\bibinfo {author} {\bibfnamefont {C.}~\bibnamefont
  {Flynn}},\ }\href {https://pypi.org/project/fbm/} {\emph {\bibinfo {title}
  {Fractional Brownian motion realizations}}},\ \bibinfo {address} {The Python
  Package Index} (\bibinfo {year} {2019 [Online]})\BibitemShut {NoStop}%
\bibitem [{\citenamefont {Davis}\ and\ \citenamefont
  {Harte}(1987)}]{DavisHarteAlgorithm1987}%
  \BibitemOpen
  \bibfield  {author} {\bibinfo {author} {\bibfnamefont {R.~B.}\ \bibnamefont
  {Davis}}\ and\ \bibinfo {author} {\bibfnamefont {D.}~\bibnamefont {Harte}},\
  }\bibfield  {title} {\bibinfo {title} {Tests for hurst effect},\ }\href@noop
  {} {\bibfield  {journal} {\bibinfo  {journal} {Biometrika}\ }\textbf
  {\bibinfo {volume} {74}},\ \bibinfo {pages} {95} (\bibinfo {year}
  {1987})}\BibitemShut {NoStop}%
\end{thebibliography}%


\begin{thebibliography}{23}%
\makeatletter
\providecommand \@ifxundefined [1]{%
 \@ifx{#1\undefined}
}%
\providecommand \@ifnum [1]{%
 \ifnum #1\expandafter \@firstoftwo
 \else \expandafter \@secondoftwo
 \fi
}%
\providecommand \@ifx [1]{%
 \ifx #1\expandafter \@firstoftwo
 \else \expandafter \@secondoftwo
 \fi
}%
\providecommand \natexlab [1]{#1}%
\providecommand \enquote  [1]{``#1''}%
\providecommand \bibnamefont  [1]{#1}%
\providecommand \bibfnamefont [1]{#1}%
\providecommand \citenamefont [1]{#1}%
\providecommand \href@noop [0]{\@secondoftwo}%
\providecommand \href [0]{\begingroup \@sanitize@url \@href}%
\providecommand \@href[1]{\@@startlink{#1}\@@href}%
\providecommand \@@href[1]{\endgroup#1\@@endlink}%
\providecommand \@sanitize@url [0]{\catcode `\\12\catcode `\$12\catcode
  `\&12\catcode `\#12\catcode `\^12\catcode `\_12\catcode `\%12\relax}%
\providecommand \@@startlink[1]{}%
\providecommand \@@endlink[0]{}%
\providecommand \url  [0]{\begingroup\@sanitize@url \@url }%
\providecommand \@url [1]{\endgroup\@href {#1}{\urlprefix }}%
\providecommand \urlprefix  [0]{URL }%
\providecommand \Eprint [0]{\href }%
\providecommand \doibase [0]{https://doi.org/}%
\providecommand \selectlanguage [0]{\@gobble}%
\providecommand \bibinfo  [0]{\@secondoftwo}%
\providecommand \bibfield  [0]{\@secondoftwo}%
\providecommand \translation [1]{[#1]}%
\providecommand \BibitemOpen [0]{}%
\providecommand \bibitemStop [0]{}%
\providecommand \bibitemNoStop [0]{.\EOS\space}%
\providecommand \EOS [0]{\spacefactor3000\relax}%
\providecommand \BibitemShut  [1]{\csname bibitem#1\endcsname}%
\let\auto@bib@innerbib\@empty
\bibitem [{\citenamefont {Kardar}\ \emph {et~al.}(1986)\citenamefont {Kardar},
  \citenamefont {Parisi},\ and\ \citenamefont {Zhang}}]{PhysRevLett.56.889}%
  \BibitemOpen
  \bibfield  {author} {\bibinfo {author} {\bibfnamefont {M.}~\bibnamefont
  {Kardar}}, \bibinfo {author} {\bibfnamefont {G.}~\bibnamefont {Parisi}},\
  and\ \bibinfo {author} {\bibfnamefont {Y.-C.}\ \bibnamefont {Zhang}},\
  }\bibfield  {title} {\bibinfo {title} {Dynamic scaling of growing
  interfaces},\ }\href {https://doi.org/10.1103/PhysRevLett.56.889} {\bibfield
  {journal} {\bibinfo  {journal} {Phys. Rev. Lett.}\ }\textbf {\bibinfo
  {volume} {56}},\ \bibinfo {pages} {889} (\bibinfo {year} {1986})}\BibitemShut
  {NoStop}%
\bibitem [{\citenamefont {Forster}\ \emph {et~al.}(1977)\citenamefont
  {Forster}, \citenamefont {Nelson},\ and\ \citenamefont
  {Stephen}}]{PhysRevA.16.732}%
  \BibitemOpen
  \bibfield  {author} {\bibinfo {author} {\bibfnamefont {D.}~\bibnamefont
  {Forster}}, \bibinfo {author} {\bibfnamefont {D.~R.}\ \bibnamefont
  {Nelson}},\ and\ \bibinfo {author} {\bibfnamefont {M.~J.}\ \bibnamefont
  {Stephen}},\ }\bibfield  {title} {\bibinfo {title} {Large-distance and
  long-time properties of a randomly stirred fluid},\ }\href
  {https://doi.org/10.1103/PhysRevA.16.732} {\bibfield  {journal} {\bibinfo
  {journal} {Phys. Rev. A}\ }\textbf {\bibinfo {volume} {16}},\ \bibinfo
  {pages} {732} (\bibinfo {year} {1977})}\BibitemShut {NoStop}%
\bibitem [{\citenamefont {Baik}\ \emph {et~al.}(1999)\citenamefont {Baik},
  \citenamefont {Deift},\ and\ \citenamefont
  {Johansson}}]{BDJDistributionLengthLongest1999}%
  \BibitemOpen
  \bibfield  {author} {\bibinfo {author} {\bibfnamefont {J.}~\bibnamefont
  {Baik}}, \bibinfo {author} {\bibfnamefont {P.}~\bibnamefont {Deift}},\ and\
  \bibinfo {author} {\bibfnamefont {K.}~\bibnamefont {Johansson}},\ }\bibfield
  {title} {\bibinfo {title} {On the distribution of the length of the longest
  increasing subsequence of random permutations},\ }\href
  {https://doi.org/10.1090/S0894-0347-99-00307-0} {\bibfield  {journal}
  {\bibinfo  {journal} {Journal of the American Mathematical Society}\ }\textbf
  {\bibinfo {volume} {12}},\ \bibinfo {pages} {1119} (\bibinfo {year}
  {1999})}\BibitemShut {NoStop}%
\bibitem [{\citenamefont {Johansson}(2000)}]{JohShapeFluctuationsRandom2000a}%
  \BibitemOpen
  \bibfield  {author} {\bibinfo {author} {\bibfnamefont {K.}~\bibnamefont
  {Johansson}},\ }\bibfield  {title} {\bibinfo {title} {Shape {{Fluctuations}}
  and {{Random Matrices}}},\ }\href {https://doi.org/10.1007/s002200050027}
  {\bibfield  {journal} {\bibinfo  {journal} {Communications in Mathematical
  Physics}\ }\textbf {\bibinfo {volume} {209}},\ \bibinfo {pages} {437}
  (\bibinfo {year} {2000})}\BibitemShut {NoStop}%
\bibitem [{\citenamefont {Pr{\"a}hofer}\ and\ \citenamefont
  {Spohn}(2002)}]{PSScaleInvariancePNG2002}%
  \BibitemOpen
  \bibfield  {author} {\bibinfo {author} {\bibfnamefont {M.}~\bibnamefont
  {Pr{\"a}hofer}}\ and\ \bibinfo {author} {\bibfnamefont {H.}~\bibnamefont
  {Spohn}},\ }\bibfield  {title} {\bibinfo {title} {Scale invariance of the
  {{PNG}} droplet and the {{Airy}} process},\ }\href
  {https://doi.org/10.1023/A:1019791415147} {\bibfield  {journal} {\bibinfo
  {journal} {Journal of Statistical Physics}\ }\textbf {\bibinfo {volume}
  {108}},\ \bibinfo {pages} {1071} (\bibinfo {year} {2002})}\BibitemShut
  {NoStop}%
\bibitem [{\citenamefont {Cator}\ and\ \citenamefont
  {Groeneboom}(2006)}]{CatorGroen}%
  \BibitemOpen
  \bibfield  {author} {\bibinfo {author} {\bibfnamefont {E.}~\bibnamefont
  {Cator}}\ and\ \bibinfo {author} {\bibfnamefont {P.}~\bibnamefont
  {Groeneboom}},\ }\bibfield  {title} {\bibinfo {title} {{Second class
  particles and cube root asymptotics for Hammersley’s process}},\ }\href
  {https://doi.org/10.1214/009117906000000089} {\bibfield  {journal} {\bibinfo
  {journal} {The Annals of Probability}\ }\textbf {\bibinfo {volume} {34}},\
  \bibinfo {pages} {1273 } (\bibinfo {year} {2006})}\BibitemShut {NoStop}%
\bibitem [{\citenamefont {Bal\'azs}\ \emph {et~al.}(2006)\citenamefont
  {Bal\'azs}, \citenamefont {Cator},\ and\ \citenamefont
  {Sepp\"al\"ainen}}]{10.1214/EJP.v11-366}%
  \BibitemOpen
  \bibfield  {author} {\bibinfo {author} {\bibfnamefont {M.}~\bibnamefont
  {Bal\'azs}}, \bibinfo {author} {\bibfnamefont {E.}~\bibnamefont {Cator}},\
  and\ \bibinfo {author} {\bibfnamefont {T.}~\bibnamefont {Sepp\"al\"ainen}},\
  }\bibfield  {title} {\bibinfo {title} {{Cube Root Fluctuations for the Corner
  Growth Model Associated to the Exclusion Process}},\ }\href
  {https://doi.org/10.1214/EJP.v11-366} {\bibfield  {journal} {\bibinfo
  {journal} {Electronic Journal of Probability}\ }\textbf {\bibinfo {volume}
  {11}},\ \bibinfo {pages} {1094 } (\bibinfo {year} {2006})}\BibitemShut
  {NoStop}%
\bibitem [{\citenamefont {Bal{\'a}zs}\ \emph {et~al.}(2012)\citenamefont
  {Bal{\'a}zs}, \citenamefont {Komj{\'a}thy},\ and\ \citenamefont
  {Sepp{\"a}l{\"a}inen}}]{BKSMicroscopicConcavityFluctuation2012}%
  \BibitemOpen
  \bibfield  {author} {\bibinfo {author} {\bibfnamefont {M.}~\bibnamefont
  {Bal{\'a}zs}}, \bibinfo {author} {\bibfnamefont {J.}~\bibnamefont
  {Komj{\'a}thy}},\ and\ \bibinfo {author} {\bibfnamefont {T.}~\bibnamefont
  {Sepp{\"a}l{\"a}inen}},\ }\bibfield  {title} {\bibinfo {title} {Microscopic
  concavity and fluctuation bounds in a class of deposition processes},\ }\href
  {https://doi.org/10.1214/11-AIHP415} {\bibfield  {journal} {\bibinfo
  {journal} {Annales de l'Institut Henri Poincar\'e, Probabilit\'es et
  Statistiques}\ }\textbf {\bibinfo {volume} {48}},\ \bibinfo {pages} {151}
  (\bibinfo {year} {2012})}\BibitemShut {NoStop}%
\bibitem [{\citenamefont {Amir}\ \emph {et~al.}(2011)\citenamefont {Amir},
  \citenamefont {Corwin},\ and\ \citenamefont
  {Quastel}}]{ACQProbabilityDistributionFree2011a}%
  \BibitemOpen
  \bibfield  {author} {\bibinfo {author} {\bibfnamefont {G.}~\bibnamefont
  {Amir}}, \bibinfo {author} {\bibfnamefont {I.}~\bibnamefont {Corwin}},\ and\
  \bibinfo {author} {\bibfnamefont {J.}~\bibnamefont {Quastel}},\ }\bibfield
  {title} {\bibinfo {title} {Probability distribution of the free energy of the
  continuum directed random polymer in 1 + 1 dimensions},\ }\href
  {https://doi.org/10.1002/cpa.20347} {\bibfield  {journal} {\bibinfo
  {journal} {Communications on Pure and Applied Mathematics}\ }\textbf
  {\bibinfo {volume} {64}},\ \bibinfo {pages} {466} (\bibinfo {year}
  {2011})}\BibitemShut {NoStop}%
\bibitem [{\citenamefont {Matetski}\ \emph {et~al.}(2020)\citenamefont
  {Matetski}, \citenamefont {Quastel},\ and\ \citenamefont
  {Remenik}}]{MQRKPZFixedPoint2020}%
  \BibitemOpen
  \bibfield  {author} {\bibinfo {author} {\bibfnamefont {K.}~\bibnamefont
  {Matetski}}, \bibinfo {author} {\bibfnamefont {J.}~\bibnamefont {Quastel}},\
  and\ \bibinfo {author} {\bibfnamefont {D.}~\bibnamefont {Remenik}},\
  }\bibfield  {title} {\bibinfo {title} {The {{KPZ}} fixed point},\ }\href@noop
  {} {\bibfield  {journal} {\bibinfo  {journal} {arXiv:1701.00018}\ } (\bibinfo
  {year} {2020})}\BibitemShut {NoStop}%
\bibitem [{\citenamefont {Dauvergne}\ \emph {et~al.}(2018)\citenamefont
  {Dauvergne}, \citenamefont {Ortmann},\ and\ \citenamefont
  {Vir{\'a}g}}]{DOVDirectedLandscape2018}%
  \BibitemOpen
  \bibfield  {author} {\bibinfo {author} {\bibfnamefont {D.}~\bibnamefont
  {Dauvergne}}, \bibinfo {author} {\bibfnamefont {J.}~\bibnamefont {Ortmann}},\
  and\ \bibinfo {author} {\bibfnamefont {B.}~\bibnamefont {Vir{\'a}g}},\
  }\bibfield  {title} {\bibinfo {title} {The directed landscape},\ }\href@noop
  {} {\bibfield  {journal} {\bibinfo  {journal} {arXiv:1812.00309}\ } (\bibinfo
  {year} {2018})}\BibitemShut {NoStop}%
\bibitem [{\citenamefont {Bakhtin}\ and\ \citenamefont
  {Khanin}(2018)}]{BKGlobalSolutionsRandom2018}%
  \BibitemOpen
  \bibfield  {author} {\bibinfo {author} {\bibfnamefont {Y.}~\bibnamefont
  {Bakhtin}}\ and\ \bibinfo {author} {\bibfnamefont {K.}~\bibnamefont
  {Khanin}},\ }\bibfield  {title} {\bibinfo {title} {On global solutions of the
  random {{Hamilton}}\textendash{{Jacobi}} equations and the {{KPZ}} problem},\
  }\href {https://doi.org/10.1088/1361-6544/aa99a6} {\bibfield  {journal}
  {\bibinfo  {journal} {Nonlinearity}\ }\textbf {\bibinfo {volume} {31}},\
  \bibinfo {pages} {R93} (\bibinfo {year} {2018})}\BibitemShut {NoStop}%
\bibitem [{\citenamefont {Arratia}(1979)}]{ArratiaCBM}%
  \BibitemOpen
  \bibfield  {author} {\bibinfo {author} {\bibfnamefont {R.~A.}\ \bibnamefont
  {Arratia}},\ }\href@noop {} {\emph {\bibinfo {title} {Coalescing Bronwinan
  motions on the line}}}\ (\bibinfo {year} {1979})\ p.\ \bibinfo {pages}
  {134},\ \bibinfo {note} {thesis (Ph.D.)--The University of Wisconsin -
  Madison}\BibitemShut {NoStop}%
\bibitem [{\citenamefont {Charles}\ and\ \citenamefont
  {Newman}(1995)}]{CNSurfaceViewFirstPassage1995}%
  \BibitemOpen
  \bibfield  {author} {\bibinfo {author} {\bibfnamefont {M.}~\bibnamefont
  {Charles}}\ and\ \bibinfo {author} {\bibfnamefont {C.~M.}\ \bibnamefont
  {Newman}},\ }\bibfield  {title} {\bibinfo {title} {A {{Surface View}} of
  {{First}}-{{Passage Percolation}}},\ }in\ \href
  {https://doi.org/10.1007/978-3-0348-9078-6_30} {\emph {\bibinfo {booktitle}
  {Proceedings of the {{International Congress}} of {{Mathematicians}}: August
  3--11, 1994 {{Z}}\"urich, {{Switzerland}}}}}\ (\bibinfo {year} {1995})\ pp.\
  \bibinfo {pages} {1017--1023}\BibitemShut {NoStop}%
\bibitem [{\citenamefont {Seppäläinen}(2020)}]{10.1214/19-AIHP1016}%
  \BibitemOpen
  \bibfield  {author} {\bibinfo {author} {\bibfnamefont {T.}~\bibnamefont
  {Seppäläinen}},\ }\bibfield  {title} {\bibinfo {title} {{Existence,
  uniqueness and coalescence of directed planar geodesics: Proof via the
  increment-stationary growth process}},\ }\href
  {https://doi.org/10.1214/19-AIHP1016} {\bibfield  {journal} {\bibinfo
  {journal} {Annales de l'Institut Henri Poincaré, Probabilités et
  Statistiques}\ }\textbf {\bibinfo {volume} {56}},\ \bibinfo {pages} {1775 }
  (\bibinfo {year} {2020})}\BibitemShut {NoStop}%
\bibitem [{\citenamefont {Nourdin}(2012)}]{NouSelectedAspectsFractional2012}%
  \BibitemOpen
  \bibfield  {author} {\bibinfo {author} {\bibfnamefont {I.}~\bibnamefont
  {Nourdin}},\ }\href@noop {} {\emph {\bibinfo {title} {Selected Aspects of
  Fractional {{Brownian}} Motion}}},\ \bibinfo {series} {Bocconi \& \&
  {{Springer}} Series}\ No.~\bibinfo {number} {4}\ (\bibinfo {year}
  {2012})\BibitemShut {NoStop}%
\bibitem [{\citenamefont {Pimentel}(2016)}]{PimDualityCoalescenceTimes2016}%
  \BibitemOpen
  \bibfield  {author} {\bibinfo {author} {\bibfnamefont {L.~P.~R.}\
  \bibnamefont {Pimentel}},\ }\bibfield  {title} {\bibinfo {title} {Duality
  between coalescence times and exit points in last-passage percolation
  models},\ }\href {https://doi.org/10.1214/15-AOP1044} {\bibfield  {journal}
  {\bibinfo  {journal} {The Annals of Probability}\ }\textbf {\bibinfo {volume}
  {44}},\ \bibinfo {pages} {3187} (\bibinfo {year} {2016})}\BibitemShut
  {NoStop}%
\bibitem [{\citenamefont {Piterbarg}(1997)}]{piterbarg1997expansions}%
  \BibitemOpen
  \bibfield  {author} {\bibinfo {author} {\bibfnamefont {V.~V.}\ \bibnamefont
  {Piterbarg}},\ }\href@noop {} {\emph {\bibinfo {title} {Expansions and
  contractions of stochastic flows}}}\ (\bibinfo  {publisher} {University of
  Southern California},\ \bibinfo {year} {1997})\BibitemShut {NoStop}%
\bibitem [{\citenamefont {Piterbarg}\ and\ \citenamefont
  {Piterbarg}(1999)}]{PPIntermittencyTracerGradient1999a}%
  \BibitemOpen
  \bibfield  {author} {\bibinfo {author} {\bibfnamefont {L.~I.}\ \bibnamefont
  {Piterbarg}}\ and\ \bibinfo {author} {\bibfnamefont {V.~V.}\ \bibnamefont
  {Piterbarg}},\ }\bibfield  {title} {\bibinfo {title} {Intermittency of the
  {{Tracer Gradient}}},\ }\href {https://doi.org/10.1007/s002200050580}
  {\bibfield  {journal} {\bibinfo  {journal} {Communications in Mathematical
  Physics}\ }\textbf {\bibinfo {volume} {202}},\ \bibinfo {pages} {237}
  (\bibinfo {year} {1999})}\BibitemShut {NoStop}%
\bibitem [{\citenamefont {Khanin}\ and\ \citenamefont
  {Li}(2021)}]{khanin2021asymptotic}%
  \BibitemOpen
  \bibfield  {author} {\bibinfo {author} {\bibfnamefont {K.}~\bibnamefont
  {Khanin}}\ and\ \bibinfo {author} {\bibfnamefont {L.}~\bibnamefont {Li}},\
  }\bibfield  {title} {\bibinfo {title} {On asymptotic behavior of iterates of
  piecewise constant monotone maps},\ }\href@noop {} {\bibfield  {journal}
  {\bibinfo  {journal} {arxiv:2110.09731}\ } (\bibinfo {year}
  {2021})}\BibitemShut {NoStop}%
\bibitem [{SM()}]{SM}%
  \BibitemOpen
  \href@noop {} {\emph {\bibinfo {title} {See Supplemental material [url], which includes Refs.\ [22-23]}}}\BibitemShut
  {NoStop}%
\bibitem [{\citenamefont {Flynn}(line)}]{PyPI_fbm}%
  \BibitemOpen
  \bibfield  {author} {\bibinfo {author} {\bibfnamefont {C.}~\bibnamefont
  {Flynn}},\ }\href {https://pypi.org/project/fbm/} {\emph {\bibinfo {title}
  {Fractional Brownian motion realizations}}},\ \bibinfo {address} {The Python
  Package Index} (\bibinfo {year} {2019 [Online]})\BibitemShut {NoStop}%
\bibitem [{\citenamefont {Davis}\ and\ \citenamefont
  {Harte}(1987)}]{DavisHarteAlgorithm1987}%
  \BibitemOpen
  \bibfield  {author} {\bibinfo {author} {\bibfnamefont {R.~B.}\ \bibnamefont
  {Davis}}\ and\ \bibinfo {author} {\bibfnamefont {D.}~\bibnamefont {Harte}},\
  }\bibfield  {title} {\bibinfo {title} {Tests for hurst effect},\ }\href@noop
  {} {\bibfield  {journal} {\bibinfo  {journal} {Biometrika}\ }\textbf
  {\bibinfo {volume} {74}},\ \bibinfo {pages} {95} (\bibinfo {year}
  {1987})}\BibitemShut {NoStop}%
\end{thebibliography}
\providecommand{\noopsort}[1]{}\providecommand{\singleletter}[1]{#1}%

\end{document}


\preprint{APS/123-QED}

\title{Supplemental Material: Point Fields of Last Passage Percolation and Coalescing Fractional Brownian Motions }

\author{Konstantin Khanin}
\author{Liying Li}
\author{Zhanghan Yin}
\affiliation{ Department of Mathematics, University of Toronto}

\date{\today}
\maketitle
\section{Simulation}

Our documented code generating discrete cfBm and LPP, performing statistical tests, and data banks can be accessed on our GitHub repository \footnote{
  \url{https://github.com/neo-nice-to-meet-you/cfBm-KPZ/}}.

\subsection{Generation of fBM}
We generate discretized steps of cBm using the Python library \textit{fbm} \cite{PyPI_fbm}, with specifically the Davis-Harte algorithm \cite{DavisHarteAlgorithm1987}.

\subsection{Implementation }
The length of the discrete simulations are measured in the number of discrete steps, denoted by $n$. The starting positions occupy integer points $\{-k,-k+1,...,k-1,k\}$, where we choose $k=20*
\text{round}(n^{2/3})$. 

The typical deviation of a sample path grows like $n^{2/3}$, hence the number of surviving points decay like $n^{-2/3}$. Thus, by the $n$-th step, the empirical upper and lower
fields contain approximately $40$ points.

Because both cfBM and LPP paths are homogeneous, it suffices to re-scale the point fields by $n^{2/3}$ in the end.
\subsection{Sample size}
Our statistical tests are based on a sample size of $ 500$ independent cfBM/LPP, which gives rise to approximately $500*40 = 20000$ samples for the consecutive point distance and jump-$k$ ratios.

For our experiments, cfBMs are generated up to $n=1024$ steps and LPP is generated up to $n=4096$. 

\subsection{Boundary effects}

Before we perform the K-S test, $2$ points from each end of the point fields aare removed to account for finite-size effects. For instance, one end point --- either the maximum or minimum (depending on
the definition), of the lower point fields is constant across samples, and is always equal to the maximum (resp.\ minimum) of the set of starting positions. Removing $2$ points from each end of the
point fields takes care of this problem.

\section{Full Tables}

\begin{table*}
\caption{\label{tab:table0} $p$-values of K-S tests comparing the distribution of $\delta_0$ of upper (U) and lower (L) fields between different models of cfBM.}

\begin{tabular}{|c|c|c|c|c|c|c|c|}
    \hline
    & coin-flip & $\alpha=1/2$ & $\alpha=1$& $\alpha=2$& $\alpha =10$ & $\alpha =\infty$ & regenerate  \\
    \hline
    coin-flip & & U: 0.60 &  U: 0.45 &  U: 0.74 &  U: 0.72 &  U: 0.33 &  U: 0.33 \\
    \hline
    $\alpha=1/2$ & L: 0.39 & &  U: 0.89 &  U: 0.83 &  U: 0.75 &  U: 0.65 &  U: 0.40 \\
    \hline
    $\alpha=1$ & L: 0.17 & L: 0.35 & &  U: 0.72 &  U: 0.79 &  U: 0.84 &  U: 0.13\\
    \hline
    $\alpha=2$ & L: $<$0.01 & L: 0.03& L: 0.15 & &  U: 0.99 &  U: 0.99 &  U: 0.22 \\
    \hline
    $\alpha=10$ & L: $<$0.01 & L: 0.02 & L: 0.03 & L: 0.49 & &  U: 0.98 &  U: 0.20 \\
    \hline
    $\alpha=\infty$ & L: $<$0.01 & L: 0.03 & L: 0.04 & L: 0.66 & L: 0.28 & & U: 0.12 \\
    \hline
    regenerate & L: $<$0.01 & L: $<$0.01 & L: $<$0.01 & L: $<$0.01 & L: $<$0.01 & L: $<$0.01 & \\
    \hline
\end{tabular}
\end{table*}

\begin{table*}
\caption{\label{tab:table1} $p$-values of K-S tests comparing the distribution of $r_2$ of upper (U) and lower (L) fields between different models of cfBM.}
\begin{tabular}{|c|c|c|c|c|c|c|c|}
    \hline
    & coin-flip & $\alpha=1/2$ & $\alpha=1$& $\alpha=2$& $\alpha =10$ & $\alpha =\infty$ & regenerate  \\
    \hline
    coin-flip & & U: 0.76 &  U: 0.58 &  U: 0.79 &  U: 0.56 &  U: 0.16 &  U: 0.16 \\
    \hline
    $\alpha=1/2$ & L: 0.94 & &  U: 0.40 &  U: 0.91 &  U: 0.95 &  U: 0.66 &  U: 0.55 \\
    \hline
    $\alpha=1$ & L: 0.87 & L: 0.98 & &  U: 0.50 &  U: 0.30 &  U: 0.52 &  U: 0.17 \\
    \hline
    $\alpha=2$ & L: 0.04 & L: 0.04 & L: 0.12 & &  U: 0.52 &  U: 0.52 &  U: 0.58\\
    \hline
    $\alpha=10$ & L: 0.02 & L: 0.04 & L: 0.11 & L: 0.70 & &  U: 0.38 &  U: 0.36\\
    \hline
    $\alpha=\infty$ & L: 0.04 & L: 0.07 & L: 0.09 & L: 0.93 & L: 0.79 & & U: 0.82 \\
    \hline
    Regenerate & L: $<$0.01 & L: $<$0.01 & L: $<$0.01 & L: $<$0.01 & L: $<$0.01 & L: $<$0.01 & \\
    \hline
\end{tabular}
\end{table*}

\begin{table*}
\caption{\label{tab:table4} $p$-values of K-S tests comparing the distribution of $\delta_0$ and $r_k$, $1\le k\le 6$, between the upper and lower fields of every processes, in order to test symmetry.}
\begin{tabular}{|c|c|c|c|c|c|c|c|c|c|}
    \hline
    & coin-flip & $\alpha=0.5$ & $\alpha=1$ & $\alpha=2$ & $\alpha=10$ & $\alpha=\infty$ & regenerate & LPP\\
    \hline
    Distance & 0.75 & 0.15 & 0.16 & $<$0.01 & $<$0.01 & $<$0.01 & 0.02 & 0.92\\
    \hline
    Jump-1 & 0.51 & 0.17 & 0.67 & $<$0.01 & $<$0.01 & $<$0.01 & 0.22 & 0.66\\
    \hline
    Jump-2 & 0.25 & 0.58 & 0.87 & $<$0.01 & 0.01 & $<$0.01 & 0.12 & 0.73\\
    \hline
    Jump-3 & 0.55 & 0.17 & 0.51 & 0.01 & $<$0.01 & $<$0.01 & 0.52 & 0.99\\
    \hline
    Jump-4 & 0.96 & 0.35 & 0.44 & $<$0.01 & $<$0.01 & $<$0.01 & 0.03 & 0.95\\
    \hline
    Jump-5 & 0.99 & 0.20 & 0.65 & 0.02 & 0.02 & 0.03 & 0.21 & 0.98\\
    \hline
    Jump-6 & 0.63 & 0.80 & 0.42 & 0.02 & $<$0.01 & 0.01 & 0.11 & 0.42\\
    \hline
\end{tabular}
\end{table*}

\begin{table*}
\caption{\label{tab:table2} $p$-values of K-S tests comparing the  distribution of $\delta_0$ and $r_k$, $1\le k\le 6$, of LPP upper field against upper fields of different cfBM models.}
\begin{tabular}{|c|c|c|c|c|c|c|c|}
    \hline
    & coin-flip & $\alpha=1/2$ & $\alpha=1$ & $\alpha=2$ & $\alpha=10$ & $\alpha=\infty$ & regenerate \\
    \hline
    vs. LPP Upper Distance & 0.50 & 0.57 & 0.23 & 0.27 & 0.15 & 0.15 & 0.77 \\
    \hline
    vs. LPP Upper Jump-1 Ratio & 0.86 & 0.82 & 0.45 & 0.99 & 0.72 & 0.59 & 0.45\\
    \hline 
    vs. LPP Upper Jump-2 Ratio & 0.17 & 0.24 & 0.15 & 0.70 & 0.21 & 0.60 & 0.99\\
    \hline 
    vs. LPP Upper Jump-3 Ratio & 0.41 & 0.45 & 0.33 & 0.49 & 0.74 & 0.79 & 0.97\\
    \hline 
    vs. LPP Upper Jump-4 Ratio & 0.43 & 0.25 & 0.44 & 0.30 & 0.87 & 0.16 & 0.90\\
    \hline 
    vs. LPP Upper Jump-5 Ratio & 0.64 & 0.71 & 0.46 & 0.77 & 0.38 & 0.50 & 0.61\\
    \hline 
    vs. LPP Upper Jump-6 Ratio & 0.41 & 0.46 & 0.32 & 0.42 & 0.34 & 0.71 & 0.78\\
    \hline
\end{tabular}
\end{table*}

\begin{table*}
  \caption{\label{tab:table3}  $p$-values of K-S tests comparing the  distribution of $\delta_0$ and $r_k$, $1\le k\le 6$, of LPP lower field against lower fields of different cfBM models.}

\begin{tabular}{|c|c|c|c|c|c|c|c|}
    \hline
    & coin-flip & $\alpha=1/2$ & $\alpha=1$ & $\alpha=2$ & $\alpha=10$ & $\alpha=\infty$ & regenerate \\
    \hline
    vs. LPP Lower Distance & 0.10 & 0.05 & 0.01 & $<$0.01 & $<$0.01 & $<$0.01 & $<$0.01 \\
    \hline 
    vs. LPP Lower Jump-1 Ratio & 0.91 & 0.18 & 0.09 & $<$0.01 & $<$0.01 & $<$0.01 & $<$0.01\\
    \hline 
    vs. LPP Lower Jump-2 Ratio & 0.24 & 0.21 & 0.04 & $<$0.01 & $<$0.01 & $<$0.01 & 0.01\\
    \hline 
    vs. LPP Lower Jump-3 Ratio & 0.19 & 0.04 & 0.01 & $<$0.01 & $<$0.01 & $<$0.01 & 0.28\\
    \hline 
    vs. LPP Lower Jump-4 Ratio & 0.45 & 0.32 & 0.26 & $<$0.01 & $<$0.01 & $<$0.01 & 0.02\\
    \hline 
    vs. LPP Lower Jump-5 Ratio & 0.50 & 0.09 & 0.08 & $<$0.01 & $<$0.01 & $<$0.01 & 0.03\\
    \hline 
    vs. LPP Lower Jump-6 Ratio & 0.42 & 0.34 & 0.25 & $<$0.01 & $<$0.01 & $<$0.01 & 0.04\\
    \hline
\end{tabular}
\end{table*}

\bibliography{fBM}